\documentclass{gtart}     
\usepackage{amssymb} 
\usepackage{amsmath,amscd,epsf,eucal,latexsym} 
\newtheorem{theorem}{Theorem}[section]
\newtheorem{lemma}[theorem]{Lemma}
\newtheorem{cor}[theorem]{Corollary}
\newtheorem{prop}[theorem]{Proposition}

\theoremstyle{definition} 
\newtheorem*{rem}{Remark}
\newtheorem*{rems}{Remarks}
\newtheorem{example}{Example}
\newtheorem*{defn}{Definition} 
\newtheorem*{examp}{Example}

\DeclareMathOperator{\intr}{int}
\DeclareMathOperator{\s}{s}
\DeclareMathOperator{\rank}{rank}

\newcommand{\ra}{\rightarrow} 
\newcommand{\bd}{\partial} 
\newcommand{\BB}{\mathcal{B}} 
\newcommand{\trb}{\bd_{\pitchfork}} 
\newcommand{\tb}{\bd_{\tau}}
\newcommand{\sm}{\smallsetminus} 
\newcommand{\FF}{\mathcal{F}}

\newcommand{\GG}{\mathcal{G}}
\newcommand{\LL}{\mathcal{L}} 
\newcommand{\R}{\mathbb{R}}
\newcommand{\RR}{\mathbf{R}}
\newcommand{\Z}{\mathbb{Z}}

\renewcommand{\ss}{\subset}
\newcommand{\0}{\emptyset}
\newcommand{\Varphi}{\Phi}

\renewcommand{\o}{\circ}
\newcommand{\hra}{\hookrightarrow}
\newcommand{\xra}{\xrightarrow}
\begin{document}
%
%
%
%
%
%
%
%
%

\title{The Sutured Thurston Norm}

\authors{John Cantwell\\Lawrence Conlon}
\address{Department of Mathematics\\ St. Louis University\\ St. 
Louis, MO 
63103}\secondaddress{Department of Mathematics\\ Washington University\\ St. 
Louis, MO 
63130}
\email{CANTWELLJC@SLU.EDU\\LC@MATH.WUSTL.EDU}

\begin{abstract}
For  sutured three-manifolds $M$, there is a \emph{sutured} Thurston norm $x^{\s}$ due to M.~Scharlemann~\cite{scharle}.  Here, we show how depth one foliations of $M$ can be useful tools for computing this norm.  This uses the relation of these foliations  with fibrations of $DM$ (the double of $M$ along the manifold $R\ss\bd M$ given by the sutured structure). We also prove and use the fact that  a natural doubling map $D_{*}:H_{2}(M,\bd M)\to H_{2}(DM,\bd DM)$ is ``norm doubling'' with respect to the norms $x^{\s}$ and $x$ on $H_{2}(M,\bd M)$ and $H_{2}(DM,\bd DM)$, respectively. All of this implies significant relations between the foliation cones of~\cite{cc:cone} and the sutured norm but, in general, these relations are difficult to pin down.
\end{abstract}



%
%
\maketitlepage    
%

%

\section{Introduction} \label{intro}

If $M$ is a compact 3-manifold, Thurston~\cite{th:norm} defines a
(semi)norm $x$ on the real vector space $H_2(M,\bd M)$ (coefficients 
$\R$
will
be understood throughout), with unit ball polyhedral, and proves:

\begin{theorem}\label{Thurston}

The fibrations of $M$ over the circle that are transverse to $\bd M$
correspond up to isotopy to the rays through lattice points in the open cones over
certain top dimensional faces \emph{(}called \emph{fibered faces}\emph{)} of the unit ball of the Thurston norm.

\end{theorem}

The cones over fibered faces of the Thurston ball will be called \emph{fibration cones}.  This is slightly misleading since the classes lying in the interior of fibration cones correspond to foliations without holonomy, ``most'' of which are dense-leaved. 

Let $(M,\gamma)$ be a compact, connected, oriented, sutured
$3$--manifold~\cite{ga1}. Write
$$\bd M=\tb M\cup\trb M.$$  This notation, introduced in earlier papers of ours and in~\cite{cancon:I}, anticipates a foliation tangent to $\tb M$ and transverse to $\trb M$.  Wherever these parts of $\bd M$ meet, $M$ has a convex corner.  This notation relates to the standard sutured manifold notation as follows:
\begin{align*}
\trb M &= \gamma=A(\gamma)\cup T(\gamma),\\
\tb M &= R(\gamma)=R_{+}\cup R_{-}.
\end{align*}
Here, $A(\gamma)$ is a union of annuli and $T(\gamma)$
is a union of tori, while $R_{\pm}$ are, respectively, the outwardly and inwardly oriented portions of $R(\gamma)$. 
The choice of orientations is part of the sutured structure and each  component of $R_{-}$ is separated from a component of $R_{+}$ by annular components of $\gamma$.  Finally, each suture is a  closed curve in the interior of a component of $A(\gamma)$, parallel to and oriented with the boundary curves of this annulus.  The union of the sutures is denoted by $s$.

  We will be interested in 
taut
foliations of $M$, hence will require
that $M$ be irreducible and, as a sutured manifold, \emph{taut}.  This latter
requirement means that each component of $\tb M$ is norm-minimizing in
$H_{2} (M,\trb M)$.  In particular, if $\sigma \ss\tb M$ is an
imbedded loop bounding a disk in $M$, it also bounds a disk in $\tb
M$. 

In~\cite{cc:cone}, we proved the following 
analog of Thurston's Theorem (Theorem~\ref{Thurston}) for depth one 
foliations.

\begin{theorem}\label{cones}
Let $(M,\gamma )$ be a compact, connected, oriented, irreducible, taut, sutured
$3$--manifold. There are finitely many closed, convex,
polyhedral cones in $H_{2} (M,\bd M)$, called \emph{foliation cones}, having
disjoint interiors and such that the taut, transversely oriented, depth one foliations of 
$(M,\gamma)$  that are transverse to $\trb M$ and have the components of $\tb M$ as sole compact leaves 
correspond to the rays through integer  lattice points of $H_{2}(M,\bd M)$ in the interior of the 
foliation cones. 
\end{theorem}

\begin{rem}
Set $M_{0}=M\sm \tb M$ and remark that 
a depth one foliation as above restricts to a fibration of $M_0$ over the 
circle.  The classes in the interior of foliation cones that are not on rays through integer lattice points correspond to foliations ``almost without holonomy'' with each leaf in $M_{0}$ dense in $M$.

\end{rem}

\begin{rem}
It is known~\cite{cc:isotopy} that the  ``foliated ray''  $\left<\FF\right>$ corresponding to the depth one foliation $\FF$ determines $\FF$ up to a $C^{0}$ isotopy that is smooth in $M_{0}$.

\end{rem}

\begin{rem}
In contrast to Thurston's result, the cones in Theorem~\ref{cones} are generally not defined by a norm.  Indeed, they are not generally symmetric with respect to multiplication by $-1$.
\end{rem}

\begin{rem}
The proof of Theorem~\ref{cones} in~\cite{cc:cone} had some serious gaps.  The authors are preparing a revised version~\cite{cone:rev} of that paper that resolves these problems.
\end{rem}

There is a seminorm $x^{\s}$ for sutured manifolds, called the
\textit{sutured} Thurston norm.  This is due to
M.~Scharlemann~\cite{scharle} and, if $s=\0$, $x^{\s}$ reduces to the usual norm 
$x$.  In
this note we develop ideas relating $x^{\s}$ to the depth one foliations classified by Theorem~\ref{cones} 
and
show how this theory can be used to compute the norm.  
This makes the  
computations of the  norm, done in the examples at the end of~\cite{cc:cone}, rigorous.  In those examples, the foliation cones are unions of cones over some faces of the Thurston ball of $x^{\s}$, but this fails
in Example~\ref{interesting} of the present paper. However, even in this example,  $x^{\s}$ is closely enough related  to the foliation cones that we are able to compute the Thurston norm.

\section{Doubling}

There are three basic topics to be treated here, namely: the doubling 
map in
singular homology, the Thurston norm in sutured manifolds and their 
doubles,
and the process of inducing fibrations in the double $DM$ from 
certain depth
one foliations on $M$.

\subsection{The Doubling Map}
If $M$ is a smooth, connected, oriented, sutured manifold, we form the
double $DM$ along $\tb M$ (assumed to be nonempty).  This is defined
in complete analogy with the usual definition of the double of a
manifold along its full boundary. Thus $DM$ is an oriented manifold
formed by taking a second copy of $M$, but with opposite orientation,
and gluing the two together via the identity map on $\tb M$.  We write
$$ DM = M\cup (-M)/\!\!\sim.
$$ There is a standard way to put a smooth, oriented structure on $DM$
so that $\bd DM$ is also smooth and so that the natural reflection map
$\rho :DM\to DM$ is an orientation--reversing diffeomorphism. This map
interchanges the corresponding points of $M$ and $-M$, hence has $\tb
M$ as its set of fixed points.

Let $S\ss M$ be a smooth, properly imbedded, oriented
surface. Reversing orientations gives $-S\ss-M$.  The
double $DS=S\cup (-S)\ss DM$ can be viewed as an oriented, properly
imbedded submanifold of $DM$.  There is a technical problem that, if
$S\cap \tb M\ne\0$, smoothness of $DS$ might fail along this set. To
avoid this, one introduces a $\rho $--invariant Riemannian metric on
$DM$.  There is a $\rho $--invariant normal neighborhood $U$ of $\tb
M$ in $DM$ and an isotopy of $S$ makes $S\cap U$ saturated by the
normal fibers of $U\cap M$.  Now $DS$ is a smooth, $\rho $--invariant
subsurface of $DM$.  Of course, if $S\cap\tb M=\0$, $DS$ is the
disjoint union of $S\text{ and }-S$. Note also that $\rho |DS$ is an
orientation--reversing diffeomorphism of this surface.

A smooth triangulation of $S$ determines a smooth triangulation of
$DS$, producing singular cycles mod the boundary in $M$ and $DM$ respectively.  The
corresponding classes $[S]\in H_{2} (M,\bd M)$ and $[DS]\in H_{2}
(DM,\bd DM)$ are well defined, independently of the choice of
triangulation. We will define a canonical ``doubling'' map $$
D_{*}:H_{2} (M,\bd M)\to H_{2} (DM,\bd DM)
$$  such that $D_{*}[S]=[DS]$ and show that this map is ``norm
doubling''.

At the level of singular chains, the map $\rho|M:M\to DM $ induces a
linear map $$
\rho _{\#}:C_{\#} (M,\bd M)\to C_{\#} (DM,\tb M\cup\bd DM)
$$ commuting with the singular boundary operator $\bd_{\#}$.  Thus, we
can define $$ D_{\#} (c) = c-\rho _{\#} (c),\quad\forall\,c\in C_{\#}
(M,\bd M),
$$ noting that this also commutes with $\bd_{\#}$.  At this point,
there is a small technical problem.  The  map $D_{*}$ induced
by $D_{\#}$ takes its image in the space $H_{*} (DM,\tb M\cup\bd DM)$, whereas
we want to interpret it as a map into the space $H_{*} (DM,\bd DM)$. The crucial
property to notice is that, if the singular chain $c$ is supported in
$\tb M$, then $D_{\#} (c)=0$.

Consider the open cover $\Varphi =\{U,V\}$ of $DM$, where $U=\intr DM$
and $V$ is a normal neighborhood of $\bd DM$  with normal fibers along
$\bd (\tb M)$ lying entirely within  $\tb M$.  Let $A=\tb M\cap V$ and
note that $\bd DM$ is a deformation retract of $A\cup\bd DM$.  
By abuse of notation, we also let $\Varphi $ denote the induced open
cover on any suspace of $DM$ and we compute singular homology on $DM$
and any of its subspaces using the $\Varphi $--small singular chain
complex $C_{\#}^{\Varphi }$. That is, each singular simplex in a chain
$c\in C_{\#}^{\Varphi }$ is supported either in $U$ or in $V$.  It is
standard that the $\Varphi $--small homology $H_{*}^{\Varphi }$ is
canonically equal to the ordinary singular homology $H_{*}$, the
equality being induced by $C_{\#}^{\Varphi }\ss C_{\#}$.

If $c\in C^{\Varphi }_{\#} (\bd M)$, then, since $D_{\#}$ annihilates
all singular simplices in $\tb M$, $D_{\#} (c)$ is a chain on
$A\cup\bd DM$.  We obtain homomorphisms 
\begin{align*}
D_{\#}&: C_{\#}^{\Varphi } (M)\to C_{\#}^{\Varphi } (DM)\\
D_{\#}&: C_{\#}^{\Varphi } (\bd M)\to C_{\#}^{\Varphi } (A\cup\bd DM),
\end{align*} of chain complexes, hence a chain homomorphism $$
D_{\#}:C_{\#}^{\Varphi } (M,\bd M)\to C_{\#}^{\Varphi } (DM,A\cup\bd 
DM).
$$  This defines $$
D_{*}:H_{*} (M,\bd M)\to H_{*} (DM,A\cup\bd DM)=H_{*} (DM,\bd DM),
$$ the desired doubling map.

\begin{rem}
The above supposes that $\tb M$ meets $\trb M$.  Otherwise, $\trb M=T(\gamma)$ and the proof  that  $D_{*}:H_{*} (M,\bd M)\to H_{*} (DM,\bd DM)$ is even easier,  not requiring the use of $\Phi$-small homology.
\end{rem}

\begin{lemma}\label{DS}
If $S\ss M$ is a properly imbedded surface, then $D_{*}[S]=[DS]$.
\end{lemma}

\begin{proof}
Indeed, if $c_{S}\in Z_{2} (M,\bd M)$ is a fundamental cycle for $S$
obtained by a smooth triangulation, it is an elementary consequence of
the orientation--reversing property of $\rho :DS\to DS$ that
$c_{S}-\rho _{\#} (c_{S})\in Z_{2} (DM,\bd DM)$ is a fundamental cycle
for $DS$.
\end{proof}

Consider the inclusion map $i:M\hookrightarrow DM$ and the induced
homomorphism $$
i^{*}:H^{1} (DM)\to H^{1} (M)
$$ in real cohomology.  Using Lefschetz duality, we view this as $$
i^{*}:H_{2} (DM,\bd DM)\to H_{2} (M,\bd M).
$$

\begin{lemma}
The composition $i^{*}\o D_{*}$ is equal to the identity on $H_{2}
(M,\bd M)$.  In particular, the doubling map is injective on $H_{2}
(M,\bd M)$.
\end{lemma}

\begin{proof}
It will be enough to prove this for elements $[S]\in H_{2} (M,\bd M)$,
where $S$ is a properly imbedded, oriented surface in $M$. Indeed,
these constitute the integer lattice in $H_{2} (M,\bd M)$. By
Lemma~\ref{DS}, we must show that $i^{*}[DS]=[S]$.  The Lefschetz
dual of $[DS] $ is represented by a 1--form $\omega $ as follows. Fix
a normal neighborhood $V$ of $DS$ in $DM$.  This can be chosen so that
$V\cap\tb M$ is saturated by normal fibers, as is $V\cap\bd DM$. The
closed form $\omega $ is supported in $V$ and has integral along each
normal fiber equal to 1.  Evidently, $V\cap M$ is a normal
neighborhood of $S$ and $\omega $ restricts in $M$ to a representative
of the Lefschetz dual of $[S]$.
\end{proof}

\begin{rem}
It is easy to give a geometric definition of $$ i^{*}:H_{2} (DM,\bd
DM)\to H_{2} (M,\bd M)$$ on each element $[\Sigma ]$ of the integer
lattice. Represent this class by a properly imbedded surface $\Sigma
\ss DM$ that is transverse to $\tb M$ and note that 
$\Sigma _{+}=\Sigma \cap M$ is a properly imbedded surface in $M$.  
Then
$i^{*}[\Sigma ]=[\Sigma _{+}]$.
\end{rem}

\subsection{The Thurston Norm}\label{tn}
Roughly speaking, we define the Thurston norm in a sutured manifold by
doubling along $\trb M$, computing the Thurston norm in the doubled
manifold, and dividing by two.  This is half the norm defined by
Scharlemann in~\cite[Definition~7.4]{scharle}. 

More precisely, let $S$ be properly imbedded as usual and
connected. By a small isotopy, $\bd S$ can be assumed to be transverse
to $\bd\trb M$ and we compute $\chi _{-}^{\s} (S)$ by doubling along
$\trb M$, computing the usual $\chi _{-}$ of the doubled surface and
dividing by two. (The superscript s stands for ``sutured''.)  One can
give an intrinsic formula for this number as follows.

The 
components of $S\cap\trb M$ are circles and/or properly imbedded
arcs in annular components of $\trb M$.  These circles need not be
essential and some of the arcs might also fail to be essential in the
sense that they start and end on the same boundary component of an
annular component in $\trb M$.  We will see that these inessential
arcs and circles can be eliminated, but for the moment they are
allowed. Let $n (S)$ denote the number of arc components of $S\cap\trb
M$. Then the reader can verify that the formula for $\chi _{-}^{\s}$ 
is 
 $$ \chi _{-}^{\s} (S) = \begin{cases} -\chi (S) +\frac{1}{2}n 
(S),&\text{if this number is
positive},\\ 0,&\text{otherwise.}
\end{cases}
$$ As usual, if $S$ is not connected, one defines $\chi _{-}^{\s} (S)$ 
as
the sum of the values on each component. If $z$ is an element of the
integer lattice 
in $H_{2} (M,\bd M)$, $x^{\s} (z)$ is defined to be the minimum value 
of
$\chi _{-}^{\s} (S)$ taken over all surfaces $S\in z$.  Continuing to
follow Thurston's lead, we extend $x^{\s}$ canonically to a 
pseudonorm on
the vector space $H_{2} (M,\bd M)$ and call this the \textit{sutured}
Thurston norm. 

\begin{rems}
Instead of computing the sutured norm by doubling in $\trb M$, one can
equally well double in $\tb M$.  Again the components of $S\cap\tb M$
are properly imbedded arcs and/or circles and the number of arc
components is the same number $n (S)$.  One then notes that $2\chi
_{-}^{\s} (S)=\chi _{-} (DS)$, where $\chi _{-} (DS)$ is defined as
for the ordinary Thurston norm. 

We further remark that, by a $\chi _{-}^{\s}$--reducing homology
and/or isotopy, $S$ can be assumed to meet each annular component of
$\trb M$ only in essential arcs, each crossing the suture once, or in
essential circles, each parallel to the suture and disjoint from it.
It can be assumed also that $S$ meets each toral component only in
essential circles, although this remark is not particularly
consequential. At any rate, $n (S)$ is now just the number of times
that $\bd S$ crosses the sutures and it is elementary that this number
is even. Thus, $\chi _{-}^{\s} (S)$ is an integer, as is $x^{\s}[S]$.

\end{rems}

\begin{examp}
A decomposing disk $\Delta $ in the sense of Gabai~\cite{ga1} has 
$\chi
_{-}^{\s} (\Delta )=0 $ if it meets the sutures twice, $\chi _{-}^{\s}
(\Delta )=1$ if it meets them four times, etc.  
\end{examp}

\begin{theorem}\label{dbl}
The map $$D_{*}:H_{2} (M,\bd M)\to H_{2} (DM,\bd DM)$$ is
norm--doubling, where the sutured Thurston norm is used on the first
space and the usual Thurston norm is used on the second. Thus, if
$B$ is the Thurston ball of $M$ and $B^{*}$ that of $DM$, then $D_{*}
(B/2)= B^{*}\cap D_{*} (H_{2} (M,\bd M))$.
\end{theorem}

\begin{proof}

It is enough to prove this on elements of the integer lattice.  Let
$[S]$ be represented by a $\chi _{-}^{\s}$--minimal surface $S$.  We
have already noted that $\chi _{-} (DS)=2 \chi _{-}^{\s} (S)$, hence
it will be enough to show that $DS$ is a $\chi _{-}$--minimal
representative of $[DS]=D_{*}[S]$.  If not, let $\Sigma \in[DS]$ have
$\chi _{-} (\Sigma )<\chi _{-} (DS)$. Isotope $\Sigma $ smoothly to be
transverse to $\tb M$ and let $\Sigma _{+}=\Sigma\cap M $ and $\Sigma
_{-}=\Sigma \cap (-M)$.  If no component of $\Sigma _{\pm}$ has
positive Euler characteristic, one verifies the relation
\begin{equation}
\chi _{-} (\Sigma )=\chi _{-}^{\s} (\Sigma _{+})+\chi _{-}^{\s}
(\Sigma _{-}).\tag{$*$}
\end{equation} 
The only possible components with positive Euler characteristic are spheres or disks. In the first case, irreducibility of $M$ permits
elimination of the offensive component.  In the second, there will be
no problem if the boundary of the disk $\Delta $ meets $\tb M$ in
arcs.  Otherwise, $\bd \Delta $ is a simple closed loop either in
$\trb M$ or $\tb M$. In the first case, $\Delta $ is also a component
of $\Sigma $ in $DM$ and has zero Thurston norm. In $M$ it has zero
sutured norm, so this case also causes no problem. In the remaining
case, $\bd \Delta \ss\tb M$ and tautness of the sutured manifold
structure, together with irreducibility, yields an isotopy of $\Sigma
$ pulling the disk $\Delta $ through $\tb M$, hence eliminating it as
a component of $\Sigma _{\pm}$. Thus $(*)$ can be assumed to hold.
Interchanging the roles of $M\text{ and }-M$, if necessary, we can
then assume that $\chi _{-}^{\s} (\Sigma _{+})<\chi _{-}^{\s}
(S)$. But $$ [S]=i^{*}[DS]=i^{*}[\Sigma ]=[\Sigma _{+}],
$$ contradicting $\chi _{-}^{\s}$--minimality of $S$ in $[S]$.
\end{proof}

\subsection{Inducing Fibrations on $DM$}

In this subsection, we assume that $M$, as a sutured manifold, is not a product  $\tb M\times I$.  This insures that
$\tb M$ cannot be a fiber in a fibration of $DM$ over the circle.
We sketch some facts that are treated in greater detail
in~\cite{cc:isotopy},~\cite{cc:surg} and~\cite{cc:cone}.  

Let $\FF$ be a smooth, depth one foliation of $M$, transverse to $\trb
M$ and having the components of $\tb M$ as sole compact leaves.  A depth one leaf $L\ss M_{0}$
determines an element $\lambda (\FF)\in H^{1} (M;\Z)$ of the integer
lattice in the real cohomology space $H^{1} (M)$ via the intersection
product with loops in $M_{0}$. This class can also be represented by a
closed, nonsingular 1-form $\omega$ on $M_{0}$ that ``blows up nicely'' at
$\tb M$ (meaning that $\omega$ becomes unbounded near $\tb M$ in such a way that the 2-plane field $\ker\omega$ extends smoothly to a 2-plane field on $M$ tangent to $\tb M$). The form $\omega $ defines
$\FF|M_{0}$, hence also determines $\FF$, and its cohomology class can
be viewed as a class on $M$ via the homotopy equivalence $M_{0}\hra M$
(the natural inclusion map). For any positive constant $a$, the form
$a \omega $ also defines $\FF$, so we obtain a ``foliated ray''
$\left<\FF\right>\ss H^{1} (M)$ corresponding to $\FF$.  This ray, in turn,
determines $\FF$ up to an isotopy that is smooth in $M_{0}$ and
continuous on $M$~\cite[Theorem~1.1]{cc:isotopy}. We often think of a
foliated ray as an isotopy class of foliations.  These foliated rays
are exactly the rays meeting integer lattice points in the interiors
of the foliation cones of~\cite{cc:cone}.  

\begin{rem}
Poincar\'e duality identifies $H^{1}(M)=H_{2}(M,\bd M)$.
\end{rem}

The leaves of $\FF|M_{0}$ spiral in a well--understood way on each 
component
$F$ of $\tb M$, giving rise to a nondivisible cohomology class $$ \nu
:\pi _{1} (F)\to \Z
$$ called the \textit{juncture} of the spiral
(cf.~\cite[\S3]{cc:isotopy}). The juncture on $F$ depends only on the
class $\lambda (\FF)$~\cite[Lemma~3.1]{cc:isotopy}.  It can be
represented by a compact, properly imbedded, oriented, nonseparating
1--manifold $N\ss F$ which need not be
connected~\cite[pp.~159--160]{cc:isotopy} and each component is
assigned an integer weight.

If there is a depth one foliation $\GG$ such that $\lambda
(\GG)=-\lambda (\FF)$, we will denote $\GG$ by $-\FF$ and call this 
the
\textit{opposite} foliation to $\FF$.  Remark that
this is not the foliation defined by the form $-\omega $, even up to
isotopy, since this foliation would require that the outwardly
oriented components of $\tb M$ become inwardly oriented and vice
versa. These orientations are part of the given sutured structure on
$M$ and may not be reversed.  While, in many cases, $-\FF$ exists,
examples show that it may not. Indeed, the three vertices in
Figure~\ref{(2,2,2)cones} of Section~\ref{exs} are not foliated classes, 
but they are the
negatives of foliated classes. Of course, at the cohomology level,
$[-\omega ]=\lambda (-\FF)$. By the ideas in the proof
of~\cite[Lemma~3.1]{cc:isotopy}, the juncture for $-\FF$ can be
represented by $-N$, the manifold obtained by reversing the
orientation of $N$.  Intuitively, the foliations $\FF\text{ and }-\FF$
spin in ``opposite directions'' along $F$, appearing to be ``mirror
images'' of one another in a small normal neighborhood of $F$ in $M$.

Suppose that $\FF$ admits an opposite foliation $-\FF$.  We can
produce a taut foliation $\FF\cup-\FF$ on $DM$ by using $\FF$ in $M$
and $-\FF$ in $-M$, the components of $\tb M$ being the sole compact
leaves. Since the foliation is taut, each of the compact leaves is a
properly imbedded, incompressible surface in $ DM$.

If $F$ is one of these compact leaves, it inherits an orientation so
that it is inwardly oriented with respect to $M$ or $-M$ and outwardly
oriented with respect to the other. Thus the junctures in $F$ for the
respective foliations can be taken to be physically the same
submanifold of $F$, but with opposite orientations. It follows that
the procedure in~\cite[pp.~379--381]{cc:surg} applies, allowing us to
erase these compact leaves by deleting their ``spiral ramp''
neighborhoods and fitting the resulting foliations together, matching
convex corners of one to concave corners of the other and vice versa
(cf.~\cite[Fig.~4]{cc:surg}).  Actually, our situation is a bit more
complicated than that envisioned in~\cite{cc:surg} because our
juncture need not be connected, but essentially the same construction
goes through. In this way we erase all leaves that are components of $\tb M$.  The
resulting foliation of $DM$, denoted by $D\FF$, has only compact
leaves since the construction amputates the finitely many ends of all
leaves and joins together their compact cores.  Thus, $D\FF$ is a
fibration of $DM$ over the circle, the fibers being transverse to $\bd
DM$.  The reader should be warned that $D\FF$ is not uniquely
determined by $\FF$ and $-\FF$.  The topology of the fiber depends on
the choices of spiral ramp neighborhoods of the components $F$ of $\tb
M$.  With a little care, this construction can be carried out so that
the following is true.

\begin{lemma}\label{fibration}
If the depth one foliation $\FF$ admits an opposite foliation, then
there are associated fibrations $D\FF$ of $DM$ over the circle 
with fibers
transverse to $\bd DM$. Furthermore, there is a smooth,
one--dimensional foliation $\LL$ of $DM$, tangent to $\bd DM$ and
transverse both to $\FF\cup-\FF$ and $D\FF$.
\end{lemma}

While each component $F$ of $\tb M$ fails to be a leaf of $D\FF$,
it remains an incompressible surface in $DM$ with a special
relationship to $D\FF$.

\begin{lemma}\label{saddles}
The surface $F$ is isotopic through properly imbedded surfaces in $DM$
to a surface that has only positive saddle tangencies with $D\FF$.
\end{lemma}

\begin{proof}
The tangent bundles $\tau =\tau (\FF\cup-\FF)$ and $\tau _{0}=\tau
(D\FF)$ are both transverse to $\LL$ and transversely oriented so
that both induce the same orientation along $\LL$. It follows that
$\tau \text{ and }\tau _{0}$ are homotopic as oriented 2--plane
bundles, hence have the same (relative) Euler class $e (\tau )=e (\tau
_{0})\in H^{2} (M,\bd M)$. Thus $$\int_{F}^{}e (\tau
_{0})=\int_{F}^{}e (\tau ) =\chi (F).$$ We can assume, via a small
isotopy near $\bd DM$, that each component of $\bd F$ is either
transverse to $D\FF$ or lies in a fiber of $D\FF$.  The two
possibilities correspond, respectively, to the cases in which the
component of $\bd F$ does or does not meet the juncture for $\FF$.
Thus, Thurston's general position result~\cite[Theorem~4]{th:norm}
allows us to perform an isotopy of $F$, putting it in a position so
that all tangencies with $D\FF$ are saddles. (The possibility that
$F$ could be isotoped onto a fiber is eliminated by our assumption
that $M$ is not a product.) If some tangency is not positive (that is,
the orientations of $\tau (F)$ and $\tau (D\FF)$ at the tangency
are opposite), it would follow that $\int_{F}^{}e (\tau _{0})\ne\chi
(F)$, a contradiction.
\end{proof}

\begin{rem}
Lemma~\ref{saddles} can also be proven more directly by a Morse
theoretic argument.
\end{rem}

\begin{prop}\label{positivesaddles}
If the depth one foliation $\FF$ admits an opposite foliation and if $K\ss DM$ is a properly imbedded surface having only positive saddle tangencies with $D\FF$, then $[K]\in H_{2}(DM,\bd DM)$ lies in the cone over a fibered face of the Thurston ball and $K$ is a norm minimizing representative of $[K]$.
\end{prop}

\begin{proof}
Let $C\ss H_{2}(DM,\bd DM)$ be the cone over a top dimensional face of the Thurston ball, the interior of which contains contains the ``fibered ray'' $\left<D\FF\right>$ associated to $D\FF$ as in Theorem~\ref{Thurston}.  Let $[D\FF]\in\left<D\FF\right>\sm\{0\}$. Then, by a standard argument of Thurston~\cite{th:norm}, the fact that the tangencies are positive saddles implies that  the convex combination $t[D\FF]+(1-t)[K]\in\intr C$, $0<t\le1$.  Consequently, $[K]\in C$.  The norm $x$ is linear in $C$, coinciding there with the linear functional $-e(\tau(D\FF)):H_{2}(DM,\bd DM)\to\R$, and so 
$$
x([K])=-e(\tau(D\FF))([K])=-\chi(K).
$$  This latter equality is  due to the fact that the tangencies are positive saddles~\cite{th:norm} (see also~\cite[Lemma~10.1.13]{cancon:II}).
\end{proof}

\begin{cor}\label{saddles'}
If the depth one foliation $\FF$ admits an opposite foliation and if $F\ss DM$ is as in \emph{Lemma~\ref{saddles}}, then $[F]$ lies in the cone over a lower dimensional face of a fibered face of the Thurston ball and $F$ is norm minimizing in $[F]$. 
\end{cor}

\begin{proof}
Indeed, by Proposition~\ref{positivesaddles} and Lemma~\ref{saddles},  $F$ is norm minimizing in $[F]$ and that class lies in the cone over a fibered face.  It cannot be in the interior of that cone since $F$ is not the fiber of a fibration of $DM$.
\end{proof}\vspace{.07in}

\begin{cor}~\label{DSsaddles}
If the depth one foliation $\FF$ admits an opposite foliation and if $S\ss M$ is a properly imbedded surface such that $DS$ is smooth and has only positive saddle tangencies with $D\FF$, then $x^{\s}[S]=-\frac{1}{2}\chi(DS)$ and $S\in[S]$ realizes this minimal sutured norm.
\end{cor}

\begin{proof}
Indeed, by Proposition~\ref{positivesaddles}, $DS$ is norm minimizing in $[DS]$.  The assertion follows by Lemma~\ref{DS} and Theorem~\ref{dbl}.
\end{proof}

\section{Sutured Handlebodies}\label{suturedhandle}

\begin{lemma}\label{one}
There is a canonical decomposition
$$H_{2}(DM,\bd DM)= H_{2}(M,\bd M)\oplus {\rm ker}\ i^{*},$$ where 
$H_2(M,\bd
M)$ is imbedded as the image of $D_*$.
\end{lemma}

\begin{proof}
Since $i^*\o D_*$ is the identity on $H_2(M,\bd M)$, this is 
immediate.
\end{proof}

\begin{lemma}\label{two}

$\ker i^{*} \cong H_{2}(M,\trb M)$.

\end{lemma}

\begin{proof}

By the long exact cohomology sequence of the pair $(DM,M)$
$$H^0(DM)\xra{i^*} H^0(M) \xra{\bd^*}
 H^1(DM,M) \ra H^1(DM) \xra{i^*} H^1(M)\cdots.$$
and the 
fact that $i^*:H^0(DM)\to H^0(M)$ is an isomorphism, it follows that $\bd^*(H^0(M)) = 0$.  Thus, the kernel of $i^*:H^1(DM)\to H^1(M)$ is isomorphic to 
$H^1(DM,M)$.  By
excision and homotopy invariance, this space is  isomorphic to
$H^1(-M,\tb(-M))$.  There is no harm in dropping the minus sign and 
employing
Lefschetz duality to identify this space with $H_2(M,\trb M)$. Here,
the version of Lefschetz duality we are using is the seldom quoted one
proven in~\cite[Theorem~3.43]{hatch}.
\end{proof}

Let $M$ be a sutured handlebody of genus $n$.  We will let 
$\gamma_{i}$,
$1\le i\le m$, denote the sutures and also the homology class each 
suture
represents in $H_{1}(M)$.  Let $\{\gamma'_i\}_{i=1}^m$ denote the 
basis of
$H_1(\trb M)$ represented by these sutures.  Let $X\ss M$ be a 
bouquet of
circles $\alpha_{j}\subset M$, $1\le j\le n$, that is a deformation 
retract
of $M$.  Viewing $\alpha_{j}$ as representing a homology class in 
$H_{1}(M)$
as well as a curve, one obtains a basis $\{\alpha_j\}_{j=1}^n$ of 
$H_1(M)$.

Consider the map $$
W:H_1(\trb M)\to H_1(M)
$$ induced by the inclusion $\trb M\hra M$.

\begin{lemma}\label{five}

The vector space $ H_{2}(M,\trb M)$ is canonically imbedded in the vector space
$H_1(\trb M)$
as   $\ker W$.

\end{lemma}

\begin{proof}

This follows from the long exact sequence
$$\cdots\ra 0=H_{2}(M)\ra H_{2}(M,\trb M) \xra{\bd}
 H_{1}(\trb M) \xra{W} H_{1}(M)\cdots.$$
\end{proof}

\begin{rem}

In the above long exact sequence, the map $W$ can be represented by 
the $n\times m$ matrix
$$\mathbf{W} = 
\begin{bmatrix}
w_{11}& \cdots& w_{1m} \\
\vdots && \vdots \\
w_{n1}& \cdots& w_{nm}
\end{bmatrix}
.$$ Here, we coordinatize $H_1(\trb M)$ by the basis
$\{\gamma_i'\}_{i=1}^m$ and $H_1(M)$ by $\{\alpha_j\}_{j=1}^n$.  The 
columns of
$\mathbf{W}$ are the vectors $\gamma_i$, $1\le i\le m$. The column 
rank $r$
of this matrix is the rank of the linear map $W$ and the dimension  
of the
kernel of $W$ is $d=m-r$.

\end{rem}

\begin{theorem}\label{Rd}

$H_{2}(DM,\bd DM) \cong H_{2}(M,\bd M)\oplus \R^{d}$.

\end{theorem}

\begin{proof}

Indeed,
\begin{align*}
H_2(DM,\bd DM) &\cong H_2(M,\bd M)\oplus \ker i^* & 
\text{(Lemma~\ref{one})}\\
&\cong H_2(M,\bd M) \oplus H_2(M,\trb M) & \text{(Lemma~\ref{two})}\\
&\cong H_2(M,\bd M)\oplus\ker W& \text{(Lemma~\ref{five})}
\end{align*}
\end{proof} 

Let $c$ be the number of 
components of $\tb M=R_{+}\cup R_{-}$.

\begin{theorem}

One has $d\ge c-1$, with equality if and only if the linear map $W$ has rank $m-c+1$ if and only if the identification
in \emph{Lemma~\ref{one}} is $$H_{2}(DM,\bd DM) = H_{2}(M,\bd M)\oplus
\R^{c-1}.$$ If $d=c-1$, the factor $\R^{c-1}$ is generated by the
classes represented by any $c-1$ of the components of $R_{+}\cup
R_{-}$.

\end{theorem}

\begin{proof}

The first equivalence follows  since the rank of $W$ equals $m-d$ while the second equivalence is immediate  by Theorem~\ref{Rd}. By Lemma~\ref{one}, the
factor $\R^{c-1}$ is identified in $H_{2} (DM,\bd DM)$ as $\ker i^{*}$
and it is clear that each component $N_{i}$ of $R_{+}\cup R_{-}$
determines a homology class $\nu _{i}=[N_{i}]\in\ker i^{*}$.  Thus, it
will be sufficient to show that any $c-1$ of these classes are
linearly independent. This will also show that $d\ge c-1$.

First note that the classes determined by the components of $R_{+}$
are linearly independent, as are those determined by the components of
$R_{-}$. Indeed,  there is a loop
in $DM$ having intersection number 1 with any given component of
$R_{+}$ and intersection number 0 with all others.  The same argument
works for the components of $R_{-}$, proving that there is no
nontrivial linear relation between the classes corresponding to the
components of one of $R_{\pm}$.  

Next, choosing the indexing appropriately, let $\{\nu
_{i}=[N_{i}]\}_{i=1}^{c-1}$ be a choice of $c-1$ of the classes and
let $\nu _{c}=[N_{c}]$ be the omitted one. For definiteness, suppose
that $N_{c}$ is a component of $R_{+}$.  We consider a linear relation
$$ 0=\sum_{i=1}^{c-1}a_{i}\nu _{i} $$ and show that each $a_{i}$ is
forced to be zero.  For each component $N_{i}$ of $R_{-}$, there is an
arc in $M$ from $N_{c}$ to $N_{i}$ and this doubles to a loop in $DM$
that has intersection number $a_{i}$ with the right hand side of the
above relation. Thus, $a_{i}=0$ whenever $N_{i}$ is a component of
$R_{-}$.  The above relation, therefore, involves only terms
corresponding to components of $R_{+}$.  As already observed, there is
no such nontrivial relation. An entirely similar argument works
when $N_{c}$ is a component of $R_{-}$.
\end{proof}

\begin{cor}\label{d=1}

The linear map $W$ has rank $m-1$ if and only if the identification
in \emph{Lemma~\ref{one}} is $$H_2(DM,\bd DM)= H_2(M,\bd
M)\oplus\R.$$ In this case, the factor $\R$ is generated by
$[R_{+}]=[R_-]$ and both $R_+$ and $R_-$ are connected.

\end{cor}

Let $g$ be the genus of $R_{+}\cup R_{-}$.

\begin{theorem}

$m-c+1 + g = n$.

\end{theorem}

\begin{proof}

The disjoint union of $R_{+}$ and $R_{-}$ has genus $g$. The proof consists of sequentially pasting
together adjoining components of the disjoint union of $R_{+}$ and $R_{-}$ along a common suture. This operation either reduces the number of components by one
or adds a handle. The totality of such pastings produces a surface
homeomorphic to $\bd M$, a connected surface of genus $n$.  Since
there are $c$ components, $c-1$ of the pastings along sutures reduce
the number of components and the remaining $m-(c-1)$ pastings add
handles to give a total of $m-c+1 + g$ handles.  The assertion
follows.
\end{proof}

\section{Computing the Sutured Thurston Norm}\label{sectwo}

Our goal in this section is to state and prove a proposition that can
often be used to find top dimensional faces of the Thurston ball.  It
applies to all the examples at the end  of~\cite{cc:cone} and 
Example~\ref{interesting} of Section~\ref{exs}.  We let
$\left[{\bf a}_{1},\ldots,{\bf a}_{n}\right]$ denote the closed, convex hull
of a set of points $\{{\bf a}_{1},\ldots,{\bf a}_{n}\}$ in
$H_{2}(M,\bd M)$ or $H_{2}(DM,\bd DM)$ and we let $\left<{\bf
a}_{1},\ldots,{\bf a}_{n}\right>$ be the cone with base $\left[{\bf
a}_{1},\ldots,{\bf a}_{n}\right]$ and cone point ${\bf 0}$. 

\begin{defn}

A {\it simple disk decomposition} of $M$ is a complete disk 
decomposition 
of $M$ in which all the disks are disjoint proper disks in $M$. That is we can 
assume all the disk are there at the beginning when we do the disk 
decomposition rather than having to do the disk decomposition 
sequentially.

\end{defn}

The following  lemmas are  consequences of Gabai's procedure
of disk decomposition~\cite{ga0}.  If $D_{i}\ss M$ is a disk of a 
simple disk
decomposition, we will denote the class $[D_{i}]\in H_{2} (M,\bd M)$
by ${\bf e}_{i}$.

\begin{lemma}\label{lem}

If $\{+D_1,\ldots,+D_n\}$ is a simple disk decomposition of $M$
giving the depth one foliation $\FF$, then each $D_{i}$, $1\le i \le
n$, meets $\FF$ in positive saddles. Furthermore, $\left<{\bf
e}_{1},\ldots,{\bf e}_{n}\right>$ is a subcone of a foliation cone and $\left<\FF\right>\sm\{\bf 0\}\subset \intr\left<{\bf
e}_{1},\ldots,{\bf e}_{n}\right>$.

\end{lemma}

For a proof, see~\cite[Corollary~2.8]{cone:rev}.

\begin{lemma}\label{lemone}
If $\{+D_1,\ldots,+D_n\}$ is a simple disk decomposition of $M$ 
giving the foliation $\FF$, 
then $\{+D_1,\ldots,+D_n\}$ is a 
simple disk decomposition of $-M$ giving the foliation $\FF$.  
Each $D_{i}\ss-M$, $1\le i \le n$, 
meets $\FF$ in positive
saddles. Furthermore, $\left<{\bf e}_{1},\ldots,{\bf e}_{n}\right>$ is  a 
subcone of a
 foliation cone
of $-M$.
\end{lemma}

\begin{proof}
Each $D_{i}$, $1\le i \le n$, and $\FF$ have the opposite transverse
orientation in $-M$ as in $M$, as does $R(\gamma)$.
\end{proof}

\begin{lemma}\label{lemtwo}

If $\{-D_1,\ldots,-D_n\}$ is a simple disk 
decomposition of $M$ giving the foliation $\FF$, 
then each $-D_{i}$ 
meets $\FF$ in positive 
saddles and so the cone
$\left<-{\bf e}_{1},\ldots,-{\bf e}_{n}\right> = -\left<{\bf 
e}_{1},\ldots,{\bf 
e}_{n}\right>$
is a subcone of a foliation cone of both $M$ and $-M$. 

\end{lemma}

\begin{proof}

Apply Lemma~\ref{lem} and~\ref{lemone}
\end{proof}

In the following,  a boundary component of a properly imbedded surface $S$  is said to cross the sutures essentially if its intersections with annular components of $\trb M$ are essential arcs. Indeed, a small isotopy of $S$ removes any inessential intersections of $\bd S$ with sutures.  When $S=D$ is a disk of a disk decomposition, the term ``essentially'' is redundant by Gabai's definition of disk decomposition, but we will use it anyway for emphasis.

\begin{prop}\label{basic}

If  $\{D_1,\ldots,D_n\}$ and $\{-D_1,\ldots,-D_n\}$ are 
simple disk decompositions of $M$, then there is a fibration $D\FF$ of $DM$ over the circle such that 
the surfaces $D_{i}\cup -D_{i}$, $1\le i\le n$, and $R_+$ have only 
positive saddle tangencies 
with the fibration. Further the Thurston norm of  
$D_{*}{\bf e}_{i} = [D_{i}\cup -D_{i}]\in H_{2}(DM,\bd DM)$ is the 
number 
of 
times $\bd D_{i}$ essentially crosses the sutures minus $2$ and the sutured 
Thurston norm 
of ${\bf e}_{i}$ is half this number.

\end{prop}

\begin{proof}

The disk decomposition $\{D_1,\ldots,D_n\}$ (respectively 
$\{-D_1,\ldots,-D_n\}$) gives the subcone  $\left<{\bf e}_{1},\ldots,{\bf e}_{n}\right>$ of a foliation cone 
 of $M$ (respectively, it gives the subcone   
$\left<-{\bf e}_{1},\ldots,-{\bf e}_{n}\right>$ of a foliation cone of $-M$). 
If $\left<\FF\right>\subset \intr\left<{\bf e}_{1},\ldots,{\bf e}_{n}\right>$, then $\left<-\FF\right>\subset
\intr\left<-{\bf e}_{1},\ldots,-{\bf e}_{n}\right>$. Then by 
Lemma~\ref{fibration}, 
$\FF$ and $-\FF$ can be matched up across $\tb M$ 
to give a fibration $D\FF$.  Further Lemmas~\ref{lem} and 
~\ref{lemtwo} imply that  $D_{i}\cup -D_{i}$, $1\le i\le n$, has only 
positive saddle tangencies with $D\FF$ while Lemma~\ref{saddles} implies that (after a small isotopy of $DM$ moving $D\FF$ and all $D_{i}\cup -D_{i}$) $R_+$ has only 
positive saddle tangencies with $D\FF$.

Let $b_{i}$ be the number of times $\bd D_{i}$ essentially crosses the sutures. 
Then the surface $D_{i}\cup -D_{i}$ is a punctured sphere with $b_{i}$ boundary components 
and thus  $-\chi(D_{i}\cup -D_{i})=b_{i}-2$. Since this 
surface has only positive saddle tangencies with the fibration, Proposition~\ref{positivesaddles} implies that $x(D_{*}{\bf e}_{i})=b_{i}-2$ and Corollary~\ref{DSsaddles} implies that    $x^{\s}({\bf 
e}_{i})=x(D_{*}{\bf e}_{i})/2$.
\end{proof}

In the examples we are interested in,  the matrix $\mathbf{W}$ of 
Section~\ref{suturedhandle} has rank $m-1$ so, by 
Corollary~\ref{d=1},  $\tb M$ has one positive component $R_+$ and 
one negative component 
$R_-$ and $H_{2}(DM,\bd DM) = H_{2}(M,\bd M) \oplus \R$ where the 
$\R$ factor is generated by $\RR = [R_+] = [R_-]$, and, without loss, 
we can assume
$$D_*(H_{2}(M,\bd M)) = H_{2}(M,\bd M)\oplus\{0\}\subset H_{2}(DM,\bd 
DM).$$

In the following corollary,   the integer $m$ and the matrix $\mathbf{W}$ are as in
Section~\ref{suturedhandle}.

\begin{cor} \label{cor}

All four of the cones $$\left<D_{*}{\bf 
e}_{1},\ldots,D_{*}{\bf
e}_{n},\pm\RR\right> \text{ and } \left<D_{*}(-{\bf
e}_{1}),\ldots,D_{*}(-{\bf e}_{n}),\pm\RR\right>$$ are subcones of fibration   
cones \emph{(}full-dimensional if $\rank\mathbf{W}=m-1$\emph{)} and thus each  lies in a cone  over a  fibered face of 
the Thurston ball of
$DM$. Also, the sutured Thurston norm is linear on both the cones
$\left<{\bf e}_{1},\ldots,{\bf e}_{n}\right>$, $\left<-{\bf
e}_{1},\ldots,-{\bf e}_{n}\right>$ $\subset H_{2}(M,\bd M)$ and both of 
these  cones are full-dimensional subcones of foliation cones
and are contained in cones over top dimensional faces of 
the
Thurston ball of $M$.

\end{cor}

\begin{proof}

Since $\{D_1,\ldots,D_n\}$ and $\{-D_1,\ldots,-D_n\}$ are 
simple disk decompositions of $M$ and $-M$ respectively, Proposition~\ref{basic} gives a fibration $D\FF$ meeting
the surfaces $D_{i}\cup -D_{i}$, $1\le i\le n$, and $R_+$  in positive 
saddles. Thus,  $\left<D_{*}{\bf e}_{1},\ldots,D_{*}{\bf e}_{n},\RR\right> $ is  a subcone of a  fibration cone (Proposition~\ref{positivesaddles}). If $\rank\mathbf{W}=m-1$, this cone is full-dimensional by Corollary~\ref{d=1} and the fact that $\mathbf{R}$ is not in the image of $D_{*}$. Similarly, since $\{-D_1,\ldots,-D_n\}$ and $\{D_1,\ldots,D_n\}$ are 
simple disk decompositions of $M$ and $-M$ respectively, one sees that the cone $\left<D_{*}({\bf -e}_{1}),\ldots,D_{*}({\bf -e}_{n}),\RR\right> $ is  a subcone of a fibration cone (full-dimensional if $\rank\mathbf{W}=m-1$). One obtains the other two fibration cones because the Thurston ball and its fibered faces are symmetric under multiplication by $-1$.

We prove the second part of the corollary for $\left<{\bf 
e}_{1},\ldots,{\bf e}_{n}\right>$.
The proof for the cone $-\left<{\bf e}_{1},\ldots,{\bf e}_{n}\right>$ 
is identical.
We must show that if 
$ {\bf p} = 
u\cdot{\bf p}_{1} + v\cdot{\bf p}_{2}$, with ${\bf p}$, ${\bf 
p}_{1}$, ${\bf p}_{2}$ $\in$ $\left<{\bf e}_{1},\ldots,{\bf
e}_{n}\right>$ and $u,v\in\R$ then 
$x^{\s}({\bf p}) 
= u\cdot x^{\s}({\bf p}_{1}) + v\cdot x^{\s}({\bf p}_{2}))$. 
Suppose on the 
contrary that
$ x^{\s}({\bf p}) \ne 
u\cdot x^{\s}({\bf p}_{1}) + v\cdot x^{\s}({\bf p}_{2})$. Then, by 
Theorem~\ref{dbl}, 
$ x(D_{*}{\bf p}) \ne 
u\cdot x(D_{*}{\bf p}_{1}) + v\cdot x(D_{*}{\bf p}_{2})$. This 
contradicts the 
linearity of the Thurston norm over 
faces of the Thurston ball of $DM$.

Since the sutured Thurston norm is linear on 
$\left<{\bf e}_{1},\ldots,{\bf e}_{n}\right>$, 
this is an (obviously full-dimensional)
subcone of the cone over a fibered face of the Thurston ball. It is also a subcone of a foliation cone by Lemma~\ref{lem}.
\end{proof}

\section{Examples}\label{exs}

In many case we can figure out the Thurston ball of knot or link 
complements cut apart along the Seifert surface using the methods of 
Section~\ref{sectwo}. The methods of Example~\ref{(2,2,2)ex} can be 
used to make rigorous the computations of the Thurston norm 
in~\cite[\S7]{cc:cone}.

\begin{example}\label{(2,2,2)ex}

Let $M$ be the complement of the pretzel link $(2,2,2)$ cut apart along its Seifert 
surface as in~\cite[\S7, Example 1]{cc:cone} (see Figure~\ref{(2,2,2)}). 
\begin{figure}[t]
\begin{center}
\begin{picture}(300,150)(0,70)

\epsfxsize=330pt
\epsffile{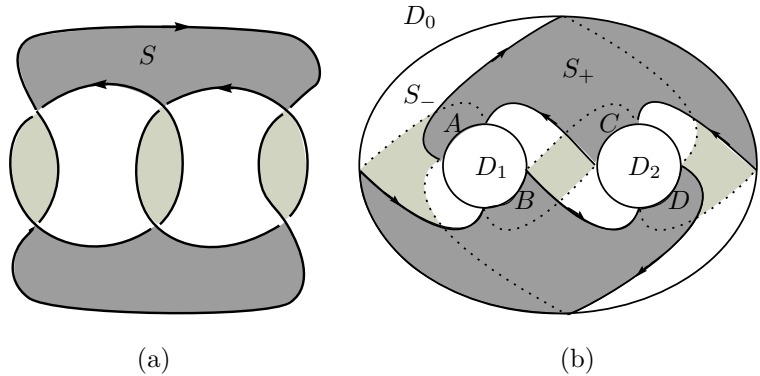}
\put (-265,70){\small (a)}
\put (-105,70){\small (b)}
\put (-265,185){\small$S$}
\put (-165,200){\small$D_{0}$}
\put (-165,170){\small$S_{-}$}
\put (-105,180){\small$S_{+}$}
\put (-150,159){\small$A$}
\put (-123,130){\small$B$}
\put (-138,143){\small$D_{1}$}
\put (-80,143){\small$D_{2}$}
\put (-91,159){\small$C$}
\put (-65,130){\small$D$}

\end{picture}
\caption{(a) A Seifert surface for $(2,2,2)$\quad (b) The sutured
manifold $M$ obtained from $(2,2,2)$}\label{(2,2,2)} 
\end{center}
\end{figure}
One can do disk 
decompositions using disks $\{D_i,-D_j\}$ as long as  $i\ne 
j\in\{0,1,2\}$. These disk decompositions are extremely easy to do 
using Gabai's graphical algorithm in~\cite[Theorem 6.1]{ga0}.  Since 
each of the $\bd D_i$'s essentially crosses the sutures $4$ times, it follows from 
Proposition~\ref{basic} that $x(D_*{\bf e}_{i})=2$ and $x^{\s}({\bf 
e}_{i})=1$.  By Corollary~\ref{cor}, it follows that the 
Thurston 
ball is the dotted hexagon $B$ of Figure~\ref{(2,2,2)cones}.  

\begin{figure}[t]
\begin{center}
\begin{picture}(300,120)(10,240)

\epsfxsize=350pt
\epsffile{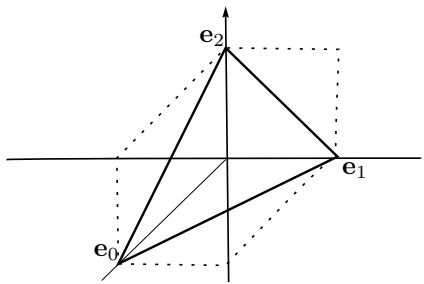}
\put (-224,259){\small$\mathbf{e}_{0}$}
\put (-130,290){\small$\mathbf{e}_{1}$}
\put (-184,340){\small$\mathbf{e}_{2}$}

\end{picture}
\caption{Thurston ball and foliation cones for 
$(2,2,2)$}\label{(2,2,2)cones}
\end{center}
\end{figure}

The Markov process 
argument of~\cite[\S7, Example~1 or Example~2]{cc:cone}  shows that 
$\left<{\bf e}_{1},{\bf 
e}_{2}\right>$, $\left<{\bf e}_{2},{\bf e}_{0}\right>$ and 
$\left<{\bf e}_{0},{\bf e}_{1}\right>$ are the foliation 
cones. 

Suitably labelling the sutures, we have that $\gamma_1 = 
-\alpha_1+\alpha_2$, $\gamma_2 = \alpha_1+\alpha_2$, and $\gamma_3 
=\alpha_1-\alpha_2$ in $H_1(M)$ (notation as in \S\ref{suturedhandle}). The matrix
$$\mathbf{W} = 
\begin{bmatrix}
-1&1&1 \\
1&1&-1
\end{bmatrix}
$$
has rank 2. 
Further, since $x(\pm\RR) =1$, Corollary~\ref{cor} allows us to conclude: \vspace{.1in}

\begin{prop}~\label{prop5-1}
The Thurston ball of $DM$ is the double cone \emph{(}suspension\emph{)} over $D_*(B/2)$ 
with cone points $\pm\RR$.
\end{prop}

\begin{rem}

Similarly, if $M$ is any of the sutured manifolds in~\cite[\S7]{cc:cone}, the 
Thur\-ston ball of $DM$ is the double cone over $D_*(B/2)$ 
with cone points $\pm \RR/x(\RR)$.

\end{rem}

\end{example}

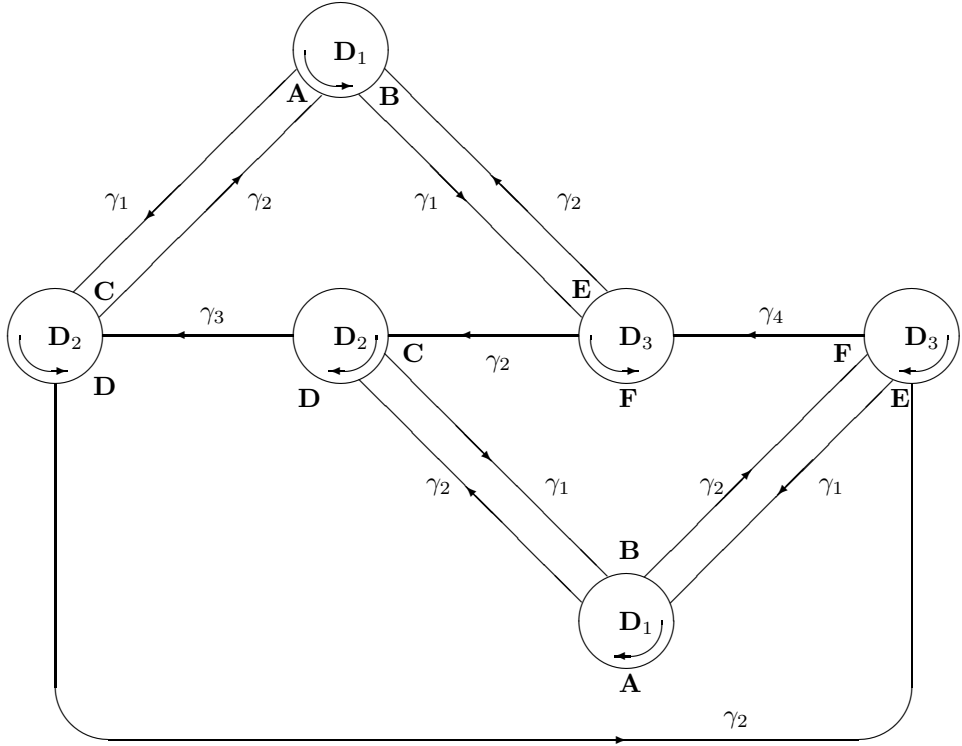
\begin{figure}[t]
\begin{center}
\setlength{\unitlength}{.9pt}
\begin{picture}(300,330)(-400,20)

\put (-434,191){\small $\mathbf D_2$}
\put(-430,195){\oval(30,30)[bl]}
\put (-430,180){\vector(1,0){5}}
\put (-430,195){\circle{40}}

\put (-314,191){\small$\mathbf D_2$}
\put(-310,195){\oval(30,30)[br]}
\put (-310,180){\vector(-1,0){5}}
\put (-310,195){\circle{40}}

\put (-194,191){\small$\mathbf D_3$}
\put(-190,195){\oval(30,30)[bl]}
\put (-190,180){\vector(1,0){5}}
\put (-190,195){\circle{40}}

\put (-74,191){\small$\mathbf D_3$}
\put(-70,195){\oval(30,30)[br]}
\put (-70,180){\vector(-1,0){5}}
\put (-70,195){\circle{40}}

\put (-314,311){\small$\mathbf D_1$}
\put(-310,315){\oval(30,30)[bl]}
\put (-310,300){\vector(1,0){5}}
\put (-310,315){\circle{40}}

\put (-194,71){\small$\mathbf D_1$}
\put(-190,75){\oval(30,30)[br]}
\put (-190,60){\vector(-1,0){5}}
\put (-190,75){\circle{40}}

\put (-334,293){\small$\mathbf {A}$}
\put (-295,292){\small$\mathbf {B}$}
\put (-415,210){\small$\mathbf {C}$}
\put (-415,170){\small$\mathbf {D}$}
\put (-214,211){\small$\mathbf {E}$}
\put (-194,165){\small$\mathbf {F}$}

\put (-194,101){\small$\mathbf{B}$}
\put (-194,45){\small$\mathbf {A}$}
\put (-329,165){\small$\mathbf {D}$}
\put (-285,184){\small$\mathbf{C}$}
\put (-80,165){\small$\mathbf {E}$}
\put (-104,184){\small$\mathbf {F}$}

\put (-422.3,213.5){\line(1,1){93.5}}
\put (-377.3,258.5){\vector(-1,-1){15}}
\put (-411.5,202.7){\line(1,1){93.5}}
\put (-366.5,247.7){\vector(1,1){15}}

\put (-291.5,307.3){\line(1,-1){93.5}}
\put (-231.5,247.3){\vector(-1,1){15}}
\put (-302.3,296.5){\line(1,-1){93.5}}
\put (-272.3,266.5){\vector(1,-1){15}}

\put (-291.5,187.3){\line(1,-1){93.5}}
\put (-261.5,157.3){\vector(1,-1){15}}
\put (-302.3,176.5){\line(1,-1){93.5}}
\put (-242.3,116.5){\vector(-1,1){15}}

\put (-171.5,82.7){\line(1,1){93.5}}
\put (-111.5,142.7){\vector(-1,-1){15}}
\put (-182.3,93.5){\line(1,1){93.5}}
\put (-152.3,123.5){\vector(1,1){15}}

\put (-90,195){\line(-1,0){80}}
\put (-120,195){\vector(-1,0){20}}
\put (-210,195){\line(-1,0){80}}
\put (-240,195){\vector(-1,0){20}}
\put (-330,195){\line(-1,0){80}}
\put (-360,195){\vector(-1,0){20}}

\put(-390,175){\oval(80,300)[bl]}
\put (-390,25){\line(1,0){280}}
\put (-210,25){\vector(1,0){20}}
\put(-110,175){\oval(80,300)[br]}


\put (-225,130){\small$\mathbf {\gamma}_1$}
\put (-110,130){\small$\mathbf {\gamma}_1$}
\put (-275,130){\small$\mathbf {\gamma}_2$}
\put (-160,130){\small$\mathbf{\gamma}_2$}

\put (-370,202){\small$\mathbf {\gamma}_3$}
\put (-250,182){\small$\mathbf{\gamma}_2$}
\put (-135,202){\small$\mathbf{\gamma}_4$}

\put (-280,250){\small$\mathbf{\gamma}_1$}
\put (-410,250){\small$\mathbf{\gamma}_1$}
\put (-350,250){\small$\mathbf{\gamma}_2$}
\put (-220,250){\small$\mathbf{\gamma}_2$}

\put (-150,32){\small$\mathbf{\gamma}_2$}
\end{picture}

\caption{A sutured handlebody}\label{example}
\end{center}
\end{figure}

\begin{example}\label{interesting}

Regard Figure~\ref{example} as drawn on $S^{2}$, the boundary of a 
solid 
ball $\mathcal{B}$. 
Paste $D_{1}$ to $D_{1}$ so that $A$ (respectively $B$) on one copy of 
$D_{1}$ 
is matched to $A$ (respectively $B$) on the other copy of $D_{1}$ 
and the sutures match up,
paste $D_{2}$ to $D_{2}$  so that $C$ (respectively $D$) on one copy of 
$D_{2}$ 
is matched to $C$ (respectively $D$) on the other copy of $D_{2}$ and the 
sutures match 
up, and paste $D_{3}$ to $D_{3}$ so that $E$ (respectively $F$) on one copy 
of $D_{3}$ 
is matched to $E$ (respectively $F$) on the other copy of $D_{3}$ 
and the sutures match up.  Then Figure~\ref{example}
represents a sutured handlebody $M$ of genus $3$ with sutures 
$\gamma_1,\gamma_2,\gamma_3,\gamma_4$. Clearly, $H_2(M,\bd M) = 
\R^3$. 

The arrows on the disks in Figure~\ref{example} define the positive 
orientation of the disks. Let $\alpha$ be a simple closed curve in 
Figure~\ref{example} going once around $D_{1}$, $D_{2}$, and $D_{3}$ 
in the negative sense and essentially crossing the sutures $\gamma_{2}$ twice and 
$\gamma_{3}$ and $\gamma_{4}$ once each. Then $\alpha$ bounds an oriented  disk 
in the solid ball $\mathcal{B}$ which we will denote $D_{0}$. In 
$H_2(M,\bd M)$, 
${\bf e}_{0} + {\bf e}_{1} + {\bf e}_{2} + {\bf e}_{3} = 0$.
\subsection{The Thurston Ball}

\begin{figure}[htb]
\begin{center}
\begin{picture}(300,220)(00,95)
\epsfxsize=300pt
\epsffile{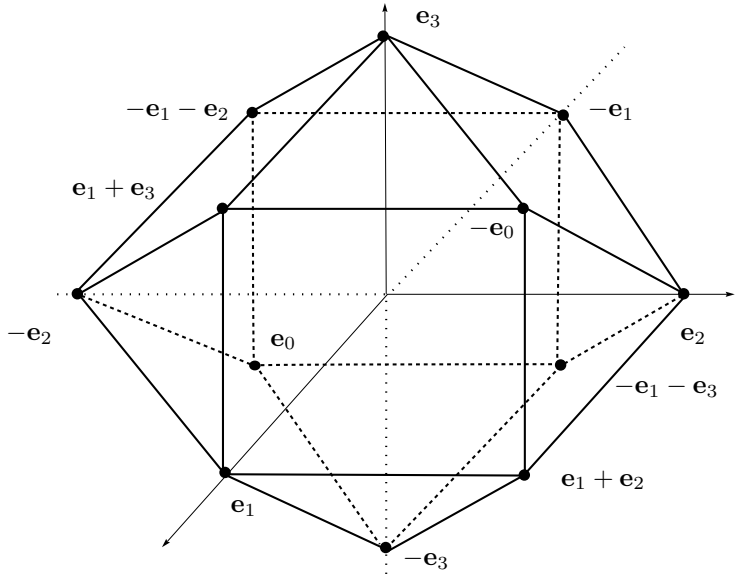}

\put (-200,183){\small$\mathbf{ e}_{0}$}
\put (-45,185){\small$\mathbf{ e}_{2}$}
\put (-145,305){\small$\mathbf {e}_{3}$}
\put (-215,120){\small$\mathbf{ e}_{1}$}

\put (-125,225){\small$-\mathbf{ e}_{0}$}
\put (-80,270){\small$-\mathbf {e}_{1}$}
\put (-300,185){\small$-\mathbf {e}_{2}$}
\put (-150,100){\small$-\mathbf{ e}_{3}$}

\put (-255,270){\small$-\mathbf{ e}_{1}-\mathbf {e}_{2}$}
\put (-275,240){\small$\mathbf{ e}_{1}+\mathbf{ e}_{3}$}
\put (-70,165){\small$-\mathbf {e}_{1}-\mathbf {e}_{3}$}

\put (-90,130){\small$\mathbf {e}_{1}+\mathbf {e}_{2}$}

\end{picture}
\caption{The Thurston ball $B$}\label{thball}
\end{center}
\end{figure}

Consider the compact, convex polyhedron depicted in Figure~\ref{thball}.  One easily checks that the vertices of the quadrilateral faces really are coplanar.  Two of these faces will present special problems in the following analysis.

\begin{defn}

The two quadrilateral faces
$Q^{\pm} = \pm[{\bf e}_{2},{\bf e}_{3},-{\bf e}_{0},-{\bf e}_{1}]$ 
 will be called the {\it exceptional faces}.

\end{defn}

\begin{lemma}
Each of the vertices in \emph{Figure~\ref{thball}} is represented by an oriented properly imbedded disk in $M$,  the boundary of which essentially crosses the sutures four times.
\end{lemma}

\begin{proof}
This is clear for $\mathbf{e}_{1}, \mathbf{e}_{2}, \mathbf{e}_{3}$ and has already been observed for $\mathbf{e}_{0}$.  For $\mathbf{e}_{1}+\mathbf{e}_{2}$, draw a closed, positively oriented curve on $\bd\BB$ meeting the suture $\gamma_{1}$ once, $\gamma_{2}$ twice, and $\gamma_{3}$ once.  This bounds the desired disk in $M$.  One argues similarly for $\mathbf{e}_{1}+\mathbf{e}_{3}$, obtaining a disk with boundary meeting $\gamma_{1}$ once, $\gamma_{2}$ twice, and $\gamma_{4}$ once.  The negatives of these classes are represented by the respective oppositely oriented disks.
\end{proof}

\begin{lemma}\label{norm}

The vertices in  \emph{Figure~\ref{thball}}  all have 
sutured Thurston norm one and the sutured norm is identically equal to $1$ on each of the nonexceptional faces.  

\end{lemma}

\begin{proof}

If $\{{\bf p}_{1}, {\bf p}_{2}, {\bf p}_{3}\}$ are any three 
vertices of a nonexceptional face with corresponding representative disks $\{\Delta_{1},\Delta_{2},\Delta_{3}\}$, then these disks and their negatives give simple disk
decompositions and each of the disks has boundary that essentially crosses  the sutures 
$4$ times.  Verifying these disk decompositions by Gabai's algorithm is routine but 
tedious. The lemma then follows by Proposition~\ref{basic} and Corollary~\ref{cor}.  
\end{proof}

Let $\pm D_{i}$ denote the disk representing $\pm\mathbf{e}_{i}$, $0\le i\le3$.
Then there are simple disk decompositions  
$\{-D_{1},D_{2},D_{3}\}$ and $\{-D_{0},D_{2},D_{3}\}$
and simple disk decompositions 
$\{D_{0},D_{1},-D_{2}\}$ and $\{D_{0},D_{1},-D_{3}\}$. There can
be no pairs of simple disk decompositions 
$\{\Delta_{1},\Delta_{2},\Delta_{3}\}$ and $\{-\Delta_{1},-\Delta_{2},-\Delta_{3}\}$ that can be 
used in Corollary~\ref{cor} to show that $Q^{\pm}$ are faces. 
Instead we will show that $x^{\s}({\bf e}_{2}+{\bf e}_{3}) =2$, which proves, by convexity of the sutured Thurston ball, that $Q^{+}$ is a face.  Of course, the norm of $-\mathbf{e}_{2}-\mathbf{e}_{3}$ is also 2 and $Q^{-}$ is a face.

In the following, $\bd (D_{2}\cup D_{3})$ and the sutures $\gamma_{i}$ are viewed as 1-cycles on $\bd M$.

\begin{lemma}\label{int}

The intersection numbers of $\bd(D_{2}\cup D_{3})$ with the sutures 
is 
given by: $\gamma_1\cdot\bd(D_2\cup D_3)=-2$, 
$\gamma_2\cdot\bd(D_2\cup D_3)=4$, 
$\gamma_3\cdot\bd(D_2\cup D_3)=-1$, $\gamma_4\cdot\bd(D_2\cup 
D_3)=-1$. 

\end{lemma}

\begin{proof}

Let ${\bf n}$ be an exterior normal to $\bd M$ and use a right hand 
rule to define the intersection number $\gamma_{i}\cdot D_{j}$, i.e. 
$\gamma_{i}\cdot D_{j} = \pm 1$ depending on whether
$(\gamma_{i},D_{j},{\bf n})$ is a right or left handed system $1\le 
i,j\le 3$.
One can compute the intersection numbers:
\begin{align*}
\gamma_{1}\cdot \bd D_{2} &= -1  & \gamma_{2}\cdot \bd D_{2} &= +2 &
\gamma_{3}\cdot \bd D_{2} &= -1  & \gamma_{4}\cdot \bd D_{2} &= 0 \\
\gamma_{1}\cdot \bd D_{3} &= -1  & \gamma_{2}\cdot \bd D_{3} &= +2 &
\gamma_{3}\cdot \bd D_{3} &= 0  & \gamma_{4}\cdot \bd D_{3} &= -1 
\end{align*}
The lemma follows.
\end{proof}

\begin{lemma}\label{disk}
If $D$ is a properly embedded disk in $M$ and $\bd D$ crosses the 
sutures essentially at most twice, 
then $D$ is boundary compressible. If $S$ 
is a properly embedded, 
connected surface in $M$ which is not a boundary compressible disk
and whose boundary 
crosses the sutures essentially \emph{(}and does so cross some sutures\emph{)}, then $ \chi_{-}(DS)\ge 2$.

\end{lemma}

\begin{proof}

Suppose $D$ is a properly embedded disk with $\bd D$ 
meeting the sutures at most twice. Put $D$ into general position with 
respect 
to $D_{1}$, $D_{2}$, and 
$D_{3}$. The points of intersection of $D$ with $D_{1}$, $D_{2}$, and 
$D_{3}$ will consist of circles and arcs. 
Assume the ends of the arcs do not lie on sutures.

By an innermost circle on $D$ 
argument, we can get rid of all circles of intersection. 

Similarly, by an innermost arc  argument on $D$ we can get rid of 
all arcs of intersection without increasing
the number of intersections of $\bd D$ 
with the sutures. In fact, choose an arc of intersection $\alpha$ in 
$D$ having endpoints $x$ and $y$
such that there exists an arc $\beta\subset \bd D$ having endpoints 
$x$ and $y$ with 
$\alpha\cup\beta$ bounding a disk $D'\subset D$ such that $\intr D'$ 
meets none of the arcs in the innermost arc argument. 
Since there are at least two such $\alpha$ and $\beta$ and since 
$\bd D$ meets the sutures at most twice, we can assume $\alpha$ 
and $\beta$ chosen so that $\beta$ meets the sutures at most once.
The arc $\alpha$ will be a properly embedded arc in 
$D_{i_{0}}$, some $1\le i_{0}\le 3$.
Thus, there is an arc 
$\delta\subset \bd D_{i_{0}}$ with endpoints $x$ and $y$, such that 
$\alpha\cup\delta$ bounds a disk $D''\subset D_{i_{0}}$. Since $\bd 
D_{i_{0}}$ meets the sutures four times and there are two possible 
choices of $\delta$, we can assume $\delta$ meets 
the sutures at most twice. Thus $\delta\cup\beta$ is a simple closed 
curve in $\bd M$ meeting the sutures at most three times, therefore never or 
twice.
Therefore $\delta\cup\beta$ bounds a disk $D'''\subset \bd M$
($D'''$ lies on the sphere represented in Figure~\ref{example} and 
$D'''$ contains none of $\pm D_{j}$, $1\le j\le 3$) and a suture 
meets $\delta$ 
if and only if it meets $\beta$.
Since $M$ is irreducible, the sphere $D'\cup D''\cup D'''$ bounds a ball that can be used to 
give an isotopy of $D$ removing the arc of intersection $\alpha$. Indeed, $D'$ can be moved onto $D''$, keeping $\alpha$ fixed, and then an arbitrarily small isotopy pulls this image of $D'$ free of $D_{i_{0}}$. Since a 
suture 
meets $\delta$ if and only if it meets $\beta$, the isotopy does not 
change the number of intersections of $\bd D$ with the sutures.
After finitely many isotopies, we may assume that $D$ does 
not meet $D_{i}$, $1\le i\le 3$ and that $\bd D$ meets the sutures 
at most twice. 
Cut $M$ apart along $D_{1}$, $D_{2}$, 
and $D_{3}$ to give the solid ball $\mathcal{B}$ with boundary 
$S^{2}$ (see 
Figure~\ref{example}).
Clearly, $D$ is boundary compressible in the solid ball 
$\mathcal{B}$ and so in $M$. 

Thus if $S$ has boundary meeting the sutures and $S$ 
is not a boundary compressible 
disk with $\bd S$ meeting the sutures twice, then either $S$ is a 
disk 
with $\bd S$ meeting the sutures 4 or more times, or 
$S$ has genus $g\ge 1$, or $S$ has at least $2$ boundary components 
and $S$ has genus $g=0$.
In the first case $\chi_{-}(DS) \ge 4-2 =2$ and, in the second case, 
$\chi_{-}(DS)\ge 2 + 4g - 2 = 4g 
> 2$.  The third case falls into two subcases.  If only one boundary component meets $\tb M$, then $DS$ has genus $0$ and at least four boundary components, in which case $ \chi_{-}(DS)\ge 4 + 0 - 2 =  2$.
If at least two boundary components of $S$ meet $\tb M$, then $DS$ has genus at least $1$ and at least two boundary components, hence $ \chi_{-}(DS)\ge 2 + 2 - 2 =  2$.
\end{proof}

\begin{lemma}\label{except}

$x^{\s}({\bf e}_{2} + {\bf e}_{3}) = 2$ and so $x^{\s}\equiv1$ on each of the exceptional faces $Q^{\pm}$.

\end{lemma}

\begin{proof}

The double of $S=D_2\cup D_3$ consists of two four times punctured 
spheres with Euler characteristic $2\cdot(2-4) = -4$. Dividing by two 
we see 
that $$x^{\s}({\bf e}_{2} + {\bf e}_{3}) \le \chi^{\s}_{-}(S) \le|-2| = 2.$$

Let $S$ be a surface representing $[D_2\cup D_3]$ in $H_{2}(M,\bd M)$.
Thus $\chi^{\s}_{-}(S) = \frac{1}{2}\chi_{-}(DS)$.
By Lemma~\ref{int},
$$\gamma_1\cdot\bd S=-2,\gamma_2\cdot\bd S=4, 
\gamma_3\cdot\bd S=-1, \gamma_4\cdot\bd S=-1.$$ 
Therefore, $\bd S$  
must meet the sutures at least eight times. 
If $S$ has only one component $S_{1}$ whose boundary meets the 
sutures,
then  $$\chi^{\s}_{-}(S) \ge \chi^{\s}_{-}(S_{1}) 
\ge \frac{1}{2}\chi_{-}(DS_{1})\ge 
\frac{1}{2}(8 + 4g -2) \ge 3,$$ where $g$ is the genus of $S_{1}$.
Otherwise $S$ has at least two 
components, $S_{1}$ and $S_{2}$, whose boundaries meet the sutures.
Thus, by Lemma~\ref{disk},  
$$\chi^{\s}_{-}(S) \ge  \frac{1}{2}\chi_{-}(DS_{1})  + 
\frac{1}{2}\chi_{-}(DS_{2}) \ge 2.$$
In any event, $x^{\s}({\bf e}_{2} + {\bf e}_{3}) \ge  2$ and equality holds. 

For the last assertion,  the fact that $x^{\s}=1$ on $\pm(\mathbf{e}_{2}+\mathbf{e}_{3})/2$  and on  each vertex of $Q^{\pm}$, together with convexity of the unit ball, implies that $x^{\s}|Q^{\pm}\equiv1$.
\end{proof}

\begin{theorem}\label{normal}

The polyhedron $B$  in 
\emph{Figure~\ref{thball}} is the unit ball of $x^{\s}$.
\end{theorem}

Indeed, by Lemma~\ref{norm} and Lemma~\ref{except}, $x^{\s}\equiv1$ on each of the faces.  

\subsection{The Foliation Cones}

\begin{figure}
\begin{center}
\begin{picture}(300,220)(00,95)
\epsfxsize=300pt
\epsffile{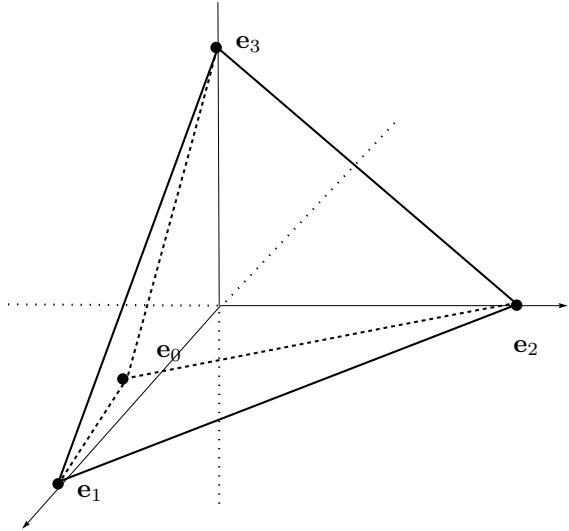}

\put (-180,183){\small$\mathbf {e}_{0}$}
\put (-45,185){\small$\mathbf {e}_{2}$}

\put (-150,300){\small$\mathbf {e}_{3}$}
\put (-210,130){\small$\mathbf {e}_{1}$}

\end{picture}
\caption{Foliation cones}\label{FolnCones}
\end{center}
\end{figure}

Bases of the foliation cones are given in  Figure~\ref{FolnCones} 
and can be found by doing the 
four simple disk decompositions using the disks 
$\{D_{1},D_{2},D_{3}\}$, 
$\{D_{0},D_{2},D_{3}\}$, $\{D_{1},D_{0},D_{3}\}$, and 
$\{D_{1},D_{2},D_{0}\}$. Thus every lattice point in the four open 
cones of Figure~\ref{FolnCones} correspond to depth one foliations.
The foliation cones obtained this way are seen to be maximal by 
the 
Markov processes argument of~\cite[\S7]{cc:cone}. 

\begin{rem}

The face $Q^{+}$ (respectively $Q^{-}$) meets the 
interior of both $\left<{\bf e}_{1},{\bf e}_{2},{\bf e}_{3}\right>$ and 
$\left<{\bf e}_{0},{\bf e}_{2},{\bf e}_{3}\right>$ (respectively
$\left<{\bf e}_{0},{\bf e}_{1},{\bf e}_{3}\right>$ and $\left<{\bf 
e}_{0},{\bf e}_{1},{\bf 
e}_{2}\right>$).
Thus none of the foliation cones can 
be the union of cones over  faces of the Thurston 
ball. 

\end{rem}

\begin{rem}
In this example it is not true that the Thurston ball of $DM$ is the double cone (suspension) of $B/2$.  The dimension of $H_{2}(DM,\bd DM)$ is 4 and $x(\pm\mathbf{R}/2)=1$ but the two exceptional faces, coned with $\pm\mathbf{R}/2$ do not give faces of the unit ball. The cones over the other faces are faces of the unit ball.
\end{rem}

\end{example}

%
%
%
%

\end{document}